\documentclass[runningheads]{llncs}
\usepackage[T1]{fontenc}
\usepackage{graphicx}
\usepackage{color}
\usepackage{graphicx}
\usepackage{epsfig}
\usepackage{tikz,float}
\usepackage{lineno}
\usepackage{amsmath,amsfonts,amssymb,mathtools}

\usepackage{geometry}

\definecolor{orcidgreen}{RGB}{166,206,57}
\newcommand{\orcidicon}{\includegraphics[width=0.26cm]{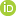}}
\def\orcidID#1{\renewcommand{\thefootnote}{\fnsymbol{footnote}}
	\unskip$^{\orcidicon}$\footnote{\orcidicon\color{orcidgreen}\,\footnotesize #1}\unskip
	\renewcommand{\thefootnote}{\arabic{footnote}}\unskip}
\begin{document}
\title{Helly Number, Radon Number and Rank in $\Delta$-Convexity on Graphs}
\author{ Bijo S Anand\inst{1}\orcidID{0000-0002-7221-904X}\and Arun Anil\inst{2}\orcidID{0000-0003-0660-5484}
 \and
Manoj Changat\inst{2}\orcidID{0000-0001-7257-6031} \and
Revathy S.Nair\inst{3}\orcidID{0009-0002-1698-5978} \and
Prasanth G. Narasimha-Shenoi\inst{4,5}\orcidID{0000-0002-5850-5410}
}
\institute{Department of Mathematics, Sree Narayana College, Punalur, 
Kollam, Kerala, India - 691305, \\ \email{bijos\_anand@yahoo.com}\and Department of Futures Studies, University of Kerala, Thiruvananthapuram, Kerala, India - 695581, \email{arunanil93@gmail.com, mchangat@gmail.com}
\and
Department of Mathematics, Mar Ivanios College, Thiruvananthapuram, Kerala, India - 695015, \email{revathyrahulnivi@gmail.com} \and
Department of Mathematics, Government College Chittur,
Palakkad, Kerala, India - 678104,  \email{prasanthgns@gmail.com}\and 
Department of Collegiate Education, Government of Kerala,
Thiruvananthapuram, Kerala, India - 695033
}
\maketitle
\begin{abstract}
This article discusses $\Delta$-convexity on simple connected graphs. We establish general bounds for the Helly number, Radon number, and rank with respect to $\Delta$-convexity on graphs. Additionally, we give the exact values for the Helly number and Radon number for chordal graphs, as well as the rank for block graphs. 

    \\
    \textbf{AMS Sub. Classification:} {52A35, 05C38, 05C69}
\end{abstract}
 \section{Introduction}\label{section_introduction}
 
	Convexity is one of the most fascinating concepts in mathematics, with connections to several fields such as optimization, combinatorics, clustering systems, and more. For further details, see \cite{boyd2004convex,rockefellar1970convex}. van de Vel \cite{van_de_vel} developed a comprehensive axiomatic treatment of convexity theory, which has  served as a foundational framework for many researchers studying convexity and its parameters.

For elaborate studies on convexity properties, refer \cite{van_de_vel}.  Several distinct notions of convexity have been introduced and investigated in the context of graphs. Pioneers such as Degreef, Doignon, Jamison, Reay, Sierksma, and van de Vel investigated different combinatorial parameters, namely, the Carath\'eodory, Helly, Radon and exchange numbers of convex sets for the axiomatic convexity; see \cite{degreef1982convex,doignon1981tverberg,jamison1974general,reay1965generalizations,sierksma1976relationships,van_de_vel}. In particular, Sierksma introduced the exchange number in 1975, \cite{sierksma1975caratheodory}.

Various authors have explored different types of convexities on graphs. For example, the geodesic convexity, see \cite{dourado2016geodetic,farber1987bridged,pelayo2013geodesic}; for monophonic convexity, see \cite{changat2004induced,duchet1988convex,morgana2002induced}; for $P_3$-convexity, see \cite{centeno2011irreversible,centeno2013geodetic}; and for triangle-path convexity, see \cite{changat2005convexities,dourado2016complexity}. 
 
A more recently studied convexity is $\Delta$-convexity, introduced in \cite{anand2020computing}, where the authors proved that computing the $\Delta$-convex hull of a general graph is NP-complete and explained various polynomial-time algorithms for some graph families like chordal, dually chordal and cographs.

The Helly number, Radon number, and Rank are fundamental concepts in convexity spaces that offer insights into their interconnectedness, separability, and dimensionality. The Helly number measures the degree of overlap among convex sets, the Radon number reflects the separability of points, and the rank indicates the maximum size of convexly independent sets (dimension of convex sets). These concepts are interrelated and have significant implications for the study of convexity spaces. Understanding them is crucial in combinatorial geometry, topology, and optimization. 
In this article, we determine the combinatorial parameters, Helly, Radon and the rank of the $\Delta$-convexity in graphs. 

Before discussing the $\Delta$-convexity, we formally define a convexity space.  A finite \emph{convexity space} is a pair $(X,\mathcal{C})$ where $X$ is a finite set and $\mathcal{C}$ is a collection of subsets of $X$ such that; $\emptyset , X \in \mathcal{C}$ and $\mathcal{C}$ is closed under intersections. The elements in $\mathcal{C}$ are said to be the \textit{convex sets}. For a subset $S$ of $X$, the smallest convex set containing $S$ is called the convex hull of $S$ and is denoted as $\langle S\rangle$.

Given a convexity $\mathcal{C}$ on $X$ and a subset $S$ of $X$, 
the set $S$ is  \emph{Helly dependent} (or, \emph{H-dependent}) provided $\displaystyle \bigcap_{a\in S}\langle S\setminus\{a\}\rangle \neq \emptyset$, and it is \emph{Helly independent} (or, \emph{H-independent})  otherwise. 
The set $S$ is \emph{Radon dependent} (or,\emph{R-dependent}) if there is a partition $\{S_1,S_2\}$ of $S$ such that $\displaystyle \langle S_1\rangle \cap \langle S_2\rangle \neq \emptyset$.
In these circumstances, $\{S_1,S_2\}$ is called a \emph{Radon partition} of $S$. If no such partition exists, then $S$ is \textit{Radon (R-) independent}.  A set $S$ is said to be \emph{Carath\'eodory dependent $\displaystyle{\langle S\rangle\subseteq \bigcup_{a\in S}\langle S\setminus\{a\}\rangle}$ and it is \emph{Carath\'eodory independent} or \textit{C-independent} otherwise}. A subset $S$ of $X$ is said to be a \textit{hull} set of $X$ if $\langle S\rangle=X$. A subset $S$ of a convex structure $X$ is \emph{convexly independent} provided $a \notin \langle S\setminus \{a\}\rangle$ for each $a \in S$. It is \emph{convexly dependent} otherwise.

A convex structure $X$ gives rise to the following numbers $h(X),r(X),c(X), d(X) \in \{0,1,\ldots, \infty\}$ determined by the following prescription. For each $n$ with $0\leq n < \infty$,
$h(X)\leq n$ if and only if each finite set $S\subseteq X$ with $|S|>n$ is H-dependent;
$r(X)\leq n$ if and only if each finite set $S\subseteq X$ with $|S|>n$ is R-dependent;
$c(X)\leq n$ if and only if each finite set $S\subseteq X$ with $|S|>n$ is C-dependent;
$d(X)\leq n$ if and only if each finite set $S\subseteq X$ with $|S|>n$ is convexly dependent;
where $h(X)$ is the \textit{Helly number} of $X$, $r(X)$ is the \textit{Radon number} of $X$, and $d(X)$ is the \emph{rank} of $X$.

In $\Delta$-convexity, the parameters Helly, Radon, Carath\'eodory and the rank are respectively denoted as $h_\Delta(X)$, $r_\Delta(X)$, $c_\Delta(X)$ and $d_\Delta(X)$ or simply $h,r,c,d$ if there is no confusion. 

 For $S\subseteq V(G)$, the \emph{$\Delta$-interval $[S]$} of $S$ is the set of all vertices in $S$ together with those vertices in $V$ that are adjacent to any two adjacent vertices in $S$.  A set $S$ is said to be \emph{$\Delta$-convex} if $[S]=S$ and the \emph{$\Delta$-convex hull} of $S$ denoted by $\langle S\rangle$, is the smallest $\Delta$-convex set containing $S$. In \cite{anand2020convexity}, the authors investigate the $\Delta$-number and $\Delta$-convexity number, providing, in particular, the exact value for the convexity number of block graphs with diameter $\leq 3$, as well as for graph products, specifically the strong and lexicographic products of two graphs, \cite{anand2022delta}.

  Throughout this article, the graphs under consideration are connected, finite, simple   and undirected, denoted by $G=(V,E)$. A graph $G$ on $n$ vertices is said to be \textit{complete} if every vertex is adjacent to every other vertex and denoted by $K_n$.  A complete graph $K_n$ is called a \textit{triangle} when $n=3$.   $G$ is said to be triangle-free if $G$ does not contain $K_3$.  A subset of vertices of a graph is \textit{independent} in a simple graph if the vertices are pairwise non-adjacent.  The \textit{independence number} $\alpha = \alpha (G)$ of a graph $G$ is the cardinality of a maximum independent set of vertices. A \textit{chordal graph} is one in which all cycles of four or more vertices have a chord, which is an edge that is not a part of the cycle but connects two non adjacent vertices of the cycle.  A \textit{block graph} is a graph in which every block is complete.  For more graph-theoretic notions that are not defined formally, see~\cite{west2001introduction}.

In this paper, for the basic definitions and results on axiomatic convexity, explained in this introductory section, we refer to \cite{van_de_vel}. The Helly number, the Radon number, and the rank of the $\Delta$-convexity are discussed in Sections~\ref{section_Helly}, \ref{section_Radon} and \ref{section_rank}, respectively. 

\begin{theorem}{\cite{van_de_vel}}\label{inequalities}
	The following holds for all convex structures.
	\begin{enumerate}
		\item (Levi inequality) $h\leq r$.
	
		\item (Eckhoff-Jamison inequality) $r \leq c(h-1)+1$ if $h\neq 1$ or $c<\infty$ 
	\end{enumerate}
\end{theorem}

\begin{remark}\cite{van_de_vel}
	A set is convexly independent if and only if all of its finite subsets are.
\end{remark} 

\begin{remark}\cite{van_de_vel} \label{rank1}
	Any convex structure satisfies $d \geq \max\{h,c,r\}$
\end{remark}
\begin{lemma}\cite{anand2020computing}\label{chordal-convex}
    If $G$ is a $2$-connected chordal graph, then every pair of adjacent vertices form a hull set of $G$.
\end{lemma}

\section{Helly number of the $\Delta$-convexity}\label{section_Helly}

In this section,  we investigate the Helly number within the context of $\Delta$-convexity. Due to the complex nature of triangle formations in graphs, determining the exact Helly number for $\Delta$-convexity presents a significant challenge. However,  we establish an upper bound for the Helly number based on the number of triangles in the graph and prove that this bound is the least upper bound. Furthermore, we explicitly derive the Helly number for chordal graphs.

We begin with the following observations that naturally arise from  the definitions of $\Delta$-convexity and Helly-independent sets:
\begin{observation}
			If $G$ is a triangle-free graph, then the Helly number, $h_\Delta(G)=|V(G)|$.
		\end{observation}
			\begin{observation}
				If $G$ is a complete graph $K_n$ for $n>2$, then $h_\Delta(G)=2$.
			\end{observation}
   The following property follows directly from the definition of a Helly dependent set: 
\begin{remark}\label{remark-helly dependent} 
	If  $S$ is a set and $w \in S$ such that $ w \in \langle S \setminus \{w\} \rangle $, then  $S$ is Helly dependent.
\end{remark}

Next, we will present some simple properties of the Helly independent set of a general graph, which will be utilized in the subsequent results.
\begin{lemma}\label{lemma_helly_characterization}
Let $G$ be a graph. Then:
\begin{enumerate}
    \item[(a)] The Helly number, $h_\Delta(G)$ satisfies $h_\Delta(G)\geq \alpha(G)$.
    \item[(b)] If $S$ is a Helly-independent set in $G$, then no three vertices of $S$ forms a $K_3$ in $G$.
    \item[(c)]If every pair of adjacent vertices forms a hull set in $G$, then $h_\Delta(G)=\max\{2,\alpha(G)\}$.
\end{enumerate}
\end{lemma}
\begin{proof}
	\begin{enumerate}
		\item[(a)] Let $S$ be a maximum independent vertex set of the graph $G$, with $|S| = \alpha(G)$. Clearly, $\langle S \rangle = S$ and $\langle S \setminus \{a\} \rangle = S \setminus \{a\}$ for all $a \in S$. This implies that $S$ is Helly-independent. Therefore, $h_\Delta(G) \geq \alpha(G)$.
		
		\item[(b)] Let $S$ be a Helly-independent set of $G$. Suppose $\{a,b,c\} \subseteq S$ forms an induced $K_3$ in $G$. Then $a \in \langle S \setminus \{a\} \rangle$, which, according to the Remark~\ref{remark-helly dependent}, implies that $S$ is Helly-dependent, a contradiction. Therefore, no three vertices of $S$ can form a $K_3$ in $G$.
 \item[(c)]	We proceed by considering the general case where  $G$ may be either a complete or a non-complete graph.
	
	First, consider the case where $G$ is a complete graph. In this case, the Helly number, $h_\Delta(G) = 2$. Since $\alpha(G) = 1$ for a complete graph, we observe that $\max\{2, \alpha(G)\} = 2$, so the lemma holds in this case.
	
	Now, suppose $G$ is not a complete graph. Let $S$ be the maximum independent vertex set of $G$. Thus $|S|=\alpha(G)$. Since every independent vertex set is H-independent, it follows that $S$ is H-independent. Now, consider any $S'\subseteq V(G)$ with $|S'|> \alpha(G)$. Since $G$ is not a complete graph, $\alpha(G)\geq 2$. Also since $|S'|\noindent > \alpha(G)$, there must exist $u,v\in S'$ with  $uv\in E(G)$ and some $x\in S'\setminus \{u,v\}$. By assumption, 
 $\langle u,v\rangle=V(G)$, implying that $\displaystyle x\in \bigcap_{a\in S'}\langle S'\setminus\{a\}\rangle$. Therefore, $S'$ is H-dependent, which implies that the Helly number, $h_\Delta(G)=\alpha(G)$.\qed
	\end{enumerate}
\end{proof}
From Lemma \ref{chordal-convex} and \ref{lemma_helly_characterization},  we immediately obtain the following corollary:

\begin{corollary}
    If $G$ is a 2-connected chordal graph, then $h_\Delta(G)=\max\{2,\alpha(G)\}$.
\end{corollary}
The following theorem provides an upper bound on the Helly number of a graph in terms of the number of triangles and the number of vertices not lying on any triangle:

\begin{theorem}
Let $G$ be a graph with $k$ triangles, and let $m$ denote the number of vertices that do not lie on any triangles. Then the Helly number, $h_\Delta(G)$ satisfies $h_\Delta(G)\leq  m + 2k$. 
\end{theorem}
\begin{proof}

Let $S \subseteq V(G)$ be a Helly independent set in $G$. By Lemma \ref{lemma_helly_characterization}, $S$ can include at most two vertices from each triangle in $G$. With \( k \) triangles in $G$, the maximum number of vertices from triangles that can be in $S$ is $2k$.
The $m$ vertices not lying on any triangle are Helly independent by default and all can be included in $S$. Therefore, the size of $S$ is bounded by $m + 2k$, implying that $h_\Delta(G) \leq m + 2k$.\qed
\end{proof}
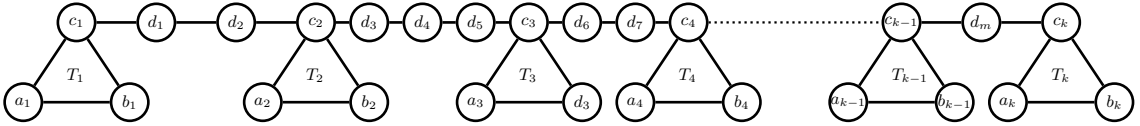
\begin{figure}[H]
		\centering
		\resizebox{1\textwidth}{!}{%
			\begin{tikzpicture}
				\tikzstyle{vertex} = [circle,fill=white,draw,minimum size=.7cm,inner sep=0pt,line width=0.05cm]
				\tikzstyle{bvertex} = [circle,fill=blue,draw,minimum size=.6cm,inner sep=0pt,thick]
				\tikzstyle{rvertex} = [circle,fill=red,draw,minimum size=.5cm,inner sep=0pt,thick]
				\tikzstyle{label} = [circle,fill=white]
				\tikzstyle{ecc} = [fill=white,thick]
				\tikzstyle{edge} = [-,line width=0.05cm]
				\tikzstyle{mvertex} = [circle,fill=magenta,draw,minimum size=.75cm,inner sep=0pt,thick]
				\tikzstyle{dvertex} = [circle,fill=white,draw,minimum size=.75cm,inner sep=0pt,thick]
                    \node[vertex] (1) at (-12,0){$c_1$};
                     \node[vertex] (2) at (-13,-1.5){$a_1$};
                     \node[vertex] (3) at (-11,-1.5){$b_1$};
                     \node[right] at (-12.3,-1) {$T_1$};
                    \node[vertex] (4) at (-10.5,0){$d_1$};
                    \node[vertex] (5) at (-9,0){$d_2$};
                     \node[vertex] (6) at (-7.5,0){$c_2$};
                      \node[vertex] (7) at (-8.5,-1.5){$a_2$};
                       \node[vertex] (8) at (-6.5,-1.5){$b_2$};
                       \node[right] at (-7.8,-1) {$T_2$};
                       \node[vertex] (9) at (-6.5,0){$d_3$};
                       \node[vertex] (10) at (-5.5,0){$d_4$};
                       \node[vertex] (11) at (-4.5,0){$d_5$};
                       \node[vertex] (12) at (-3.5,0){$c_3$};
                       \node[vertex] (13) at (-4.5,-1.5){$a_3$};
                        \node[vertex] (14) at (-2.5,-1.5){$d_3$};
                        \node[right] at (-3.8,-1) {$T_3$};
                        \node[vertex] (15) at (-2.5,0){$d_6$};
                        \node[vertex] (16) at (-1.5,0){$d_7$};
                         \node[vertex] (17) at (-.5,0){$c_4$};
                         \node[vertex] (18) at (-1.5,-1.5){$a_4$};
                         \node[vertex] (19) at (.5,-1.5){$b_4$};
                          \node[right] at (-.8,-1) {$T_4$};
                           \node[vertex] (20) at (3.5,0){$c_{k-1}$};
                             \node[vertex] (21) at (2.5,-1.5){$a_{k-1}$};
                              \node[vertex] (22) at (4.5,-1.5){$b_{k-1}$};
                              \node[right] at (3.2,-1) {$T_{k-1}$};
                              \node[vertex] (23) at (5,0){$d_{m}$};
                              \node[vertex] (24) at (6.5,0){$c_{k}$};
                              \node[vertex] (25) at (5.5,-1.5){$a_{k}$};
                              \node[vertex] (26) at (7.5,-1.5){$b_{k}$};
                               \node[right] at (6.2,-1) {$T_k$};
                              \draw[edge](1)--(2)--(3)--(1)--(4)--(5)--(6)--(7)--(8)--(6)--(9)--(10)--(11)--(12)--(13)--(14)--(12)--(15)--(16)--(17)--(18)--(19)--(17);
                              \draw[dotted,very thick](17)--(20);
                              \draw[edge](20)--(21)--(22)--(20)--(23)--(24)--(25)--(26)--(24);
			
			\end{tikzpicture}
		}
	\caption{Graph $G$ with  $d=h=r=m+2k$ }\label{fig_rank}
	
	\end{figure}

\begin{remark}
		There exists a graph $G$ with $k$ triangles  such that the Helly number,  $h_\Delta(G)=  m + 2k$, where $m$ denotes the number of vertices not lying on any triangle. For example, in Figure \ref{fig_rank}, $S=\{a_1,a_2,\ldots,a_k,b_1,b_2,\ldots,b_k,d_1,d_2,\ldots,d_m\}$, where $|S|=m+2k$. Then the following holds: 
      
      $\left\langle S\setminus\{a_i\}\right\rangle=V(G)\setminus\{a_i,c_i\}$, $\left\langle S\setminus\{b_i\}\right\rangle=V(G)\setminus\{b_i,c_i\}$ and $\left\langle S\setminus\{d_j\}\right\rangle=V(G)\setminus\{d_j\}$ for $i=1,2,\ldots,k$ and $j=1,2,3,\ldots,m$. That is, $\displaystyle{\bigcap_{a\in S}\left\langle S\setminus\{a\}\right\rangle=\emptyset}$ implying that $S$ is Helly independent. For any $S'\subseteq V(G)$ with $|S'|>m+2k$, $S'$ contains all vertices of at least one triangle. Clearly $S'$ is Helly dependent. Therefore, $h_\Delta=m+2k$. 
	\end{remark}
 
Furthermore, the Helly number can be arbitrarily large, in the sense that for any natural number $n$, one can construct a graph $G$ with $h_\Delta(G)=n$, as stated in the following proposition.
    \begin{proposition}\label{prop:hdelta}
		For any natural number $n>1$, there exists a graph $G$ with $h_\Delta(G)= n$.
	\end{proposition}
    
    \begin{proof}
 For $n=2$, consider any complete graph $G$. In this case, $h_\Delta(G)=2$.
		
		For $n\geq 3$, consider the graph $G$ depicted in Figure \ref{fig_delta}. Let $S=\{a_1,a_2,a_3,\ldots,a_n\}$. Then, $\langle S\rangle=V(G)$, $\langle S\setminus\{a_1\}\rangle=\{a_2,a_3,a_4,\ldots, a_n\}$, $\langle S\setminus\{a_2\}\rangle=\{a_1,a_3,a_4,\ldots ,a_n\}$ and $\langle S\setminus\{a_i\}\rangle=\{a_1,a_2,b_1,a_3,b_2,a_4,\ldots,a_{i-1},b_{i-2},a_{i+1},a_{i+2},\ldots ,a_n\}$ for $i=3,4,5,\ldots ,n$.
		Since $ \displaystyle \bigcap_{a\in S}\langle S\setminus\{a\}\rangle=\emptyset$, $S$ is H-independent.
		 
		 Next, we need to prove that there is no H-independent set of size greater than $n$. Suppose there exists an H-independent set $S'$ with $|S'|\geq n+1$. It is clear that $S'$ does not contains all vertices of a triangle.
   
		Now, consider the possible cases:
		 \begin{itemize}
		 	\item[(a)]Suppose $S'$ contains one vertex from each triangle. Since $G$ contains only $n-1$ triangles, $S'$ must contain two vertices from at least two triangles, say $T_i$ and $T_j$. Let  $a_i,b_i,a_j,b_j \in S'$, where $a_i,b_i \in T_i$ and $a_j,b_j \in T_j$. Then $\langle S'\setminus{a_i}\rangle=\langle S'\setminus{b_i}\rangle=\langle S'\setminus{a_j}\rangle=\langle S'\setminus{b_j}\rangle=V(G)$. By Remark \ref{remark-helly dependent},  $S'$ is H-dependent. 
		 	\item[(b)]Suppose $S'$ does not contains vertices from some triangles. In this case, $S'$ contains two vertices from at least three triangles. Therefore, there exists a sequence of triangles $T_i,T_{i+1},\ldots, T_{i+r}$ (with $i\in \{1,2,\ldots, k\}$ and $i+r\leq k$) such that $S'$ contains at least one vertex from each of these triangles and two vertices from at least two triangles, say $T_x$ and $T_y$, (with $x,y\in\{i,i+1,\ldots,i+r\}$). Let $a_x,b_x \in T_x$ and $a_y,b_y \in T_y$. Then $a_x\in\langle S'\setminus{a_x}\rangle$, $b_x\in\langle S'\setminus{b_x}\rangle$ , $a_y\in \langle S'\setminus{a_y}\rangle$ and $b_y\in \langle S'\setminus{b_y}\rangle$. This
      makes $S'$ H-dependent by Remark \ref{remark-helly dependent}. 
		 \end{itemize}
		Hence, for the graph in Figure~\ref{fig_delta},  $h_\Delta(G)=n$. \qed 

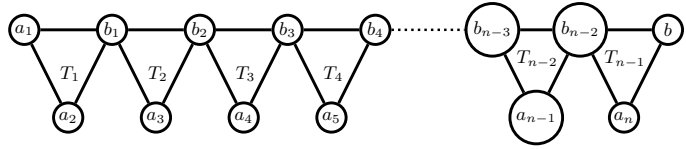
\begin{figure}[H]
		\centering
		\resizebox{0.6\textwidth}{!}{%
			\begin{tikzpicture}
				\tikzstyle{vertex} = [circle,fill=white,draw,minimum size=.5cm,inner sep=0pt,line width=0.05cm]
                \tikzstyle{lvertex} = [circle,fill=white,draw,minimum size=.9cm,inner sep=0pt,line width=0.05cm]
				\tikzstyle{bvertex} = [circle,fill=blue,draw,minimum size=.5cm,inner sep=0pt,thick]
				\tikzstyle{rvertex} = [circle,fill=red,draw,minimum size=.5cm,inner sep=0pt,thick]
				\tikzstyle{label} = [circle,fill=white]
				\tikzstyle{ecc} = [fill=white,thick]
				
				\tikzstyle{edge} = [-, line width=0.05cm]
				\tikzstyle{mvertex} = [circle,fill=magenta,draw,minimum size=.75cm,inner sep=0pt,thick]
				\tikzstyle{dvertex} = [circle,fill=white,draw,minimum size=.75cm,inner sep=0pt,thick]
				
                    \node[vertex] (1) at (-12,0){$a_1$};
                     \node[vertex] (2) at (-11.25,-1.5){$a_2$};
                 
                     \node[right] at (-11.5,-.75) {$T_1$};
                    \node[vertex] (3) at (-10.5,0){$b_1$};
                    \node[vertex] (4) at (-9,0){$b_2$};
                    \node[vertex] (5) at (-9.75,-1.5){$a_3$};
                    \node[right] at (-10,-.75) {$T_2$};
                    \node[vertex] (6) at (-7.5,0){$b_3$};
                    \node[vertex] (7) at (-8.25,-1.5){$a_4$};
                    \node[right] at (-8.5,-.75) {$T_3$};
                     \node[vertex] (8) at (-6,0){$b_4$};
                    \node[vertex] (9) at (-6.75,-1.5){$a_5$};
                    \node[right] at (-7,-.75) {$T_4$};
                    \draw[edge](1)--(2)--(3)--(1);
                    \draw[edge](3)--(4)--(5)--(3);
                    \draw[edge](4)--(6)--(7)--(4);
                    \draw[edge](6)--(8)--(9)--(6);
                     \node[lvertex] (10) at (-4,0){$b_{n-3}$};
                      \node[lvertex] (11) at (-2.5,0){$b_{n-2}$};
                      \node[right] at (-3.7,-.5) {$T_{n-2}$};
                      \node[lvertex] (12) at (-3.25,-1.5){$a_{n-1}$};
                      \node[vertex] (13) at (-1.75,-1.5){$a_{n}$};
                      \node[right] at (-2.2,-.5) {$T_{n-1}$};
                       \node[vertex] (14) at (-1,0){$b$};
                    \draw[dotted,very thick](8)--(10);
                     \draw[edge](10)--(11)--(12)--(10);
                      \draw[edge](11)--(13)--(14)--(11);
                     
			\end{tikzpicture}
		}
	\caption{Graph $G$ with with $h_\Delta(G)=d_\Delta(G)=r_\Delta(G)=n$.}\label{fig_delta}
	
	\end{figure}
\end{proof}

We conclude this section by presenting the exact value for the Helly number for chordal graphs. For the following result, we define a block of a graph $G$ as a maximal 2-connected subgraph.
\begin{theorem}
    Let $G$ be a chordal graph. Then: 
     \begin{enumerate}
        \item[(a)] If $G$  is  a block graph with $\ell$ blocks, then $h_\Delta(G) = \ell + 1$. 
      \item[(b)]  If $G$  is  not a block graph, $i.e.$, $G$ contains non-complete blocks. Let $B_1,B_2,\dots, B_\ell$ be the complete blocks and $B'_1,B'_2,\dots, B'_k$ be the non-complete blocks. Let $G'$ be the induced subgraph of $G$ formed by $\bigcup_{i=i}^{k} V(B'_j)$. Then $h_\Delta(G) =\alpha(G')+\ell $.
    \end{enumerate}
\end{theorem}
\begin{proof}
    \begin{enumerate}

        \item[(a)] Let $G$  be a chordal graph where all blocks are complete. i.e., $G$ is a block graph. Let $B_1,B_2,\dots, B_\ell$ be the blocks of $G$, with each $B_i$ being complete for $i=1,2,\dots,\ell$. Define the set $S=\{u_1,{u_1}',u_2, u_3,\ldots,u_\ell\}$, where $u_1,{u_1}'\in B_1$, $u_i \in B_i$ for $i=2,3,\dots,\ell$, and no three vertices of $S$ forms a $K_3$( a complete graph on three vertices) in $G$. Then we have the following properties:  $V(B_1)\setminus\{{u_1}'\}\nsubseteq \langle S\setminus\{u_1\}\rangle$, $V(B_1)\setminus\{{u_1}\}\nsubseteq \langle S\setminus\{u_1'\}\rangle$ and $V(B_i)\nsubseteq \langle S\setminus\{u_i\}\rangle$ for $i=2,3,\ldots,\ell$. This implies  that $\displaystyle \bigcap_{a\in S}\langle S\setminus\{a\}\rangle=\emptyset$. Hence, $S$ is H-independent.

			Next, we must prove that there is no H-independent set of size greater than $\ell+1$. Suppose,  for contradiction, that there exists an H-independent set $S'$ with $|S'| \geq \ell+2$. By Lemma \ref{lemma_helly_characterization}, $S'$ contains two vertices from at least two blocks and not more than two vertices from a single block.\\
			Consider the following possible cases:
			
			    \textbf{Case 1:} $S'$ contains at least one vertex from each block.

        Since $|S'|\geq \ell+2$, $S'$ contains two vertices from at least two blocks, say  $B_i$ and $B_j$ with $i, j \in \{1, 2, \dots, \ell\}$. Let $u, u' \in V(B_i)$ and $v, v' \in V(B_j)$, where $u, u', v, v' \in S'$. Then, $\langle S'\setminus\{u\}\rangle =\langle S'\setminus\{u'\}\rangle= \langle S'\setminus\{v\}\rangle =\langle S'\setminus\{v'\}\rangle = V(G)$. Thus, $\{u, u', v, v'\}\subseteq \langle S'\setminus\{a\}\rangle$ for all $a\in S'$, implying that $S'$ is $H$-dependent.\\
\textbf{Case 2:} $S'$ does not contains vertices from some blocks. 

   Since $|S'| \geq \ell + 2$ and $S'$ does not contains  more than two vertices from any single block, let $p$ be the number of blocks from which exactly two vertices are in $S'$ and $q$ be the number of blocks from which no vertices are in $S'$. Clearly, $p > q$.
   Thus, there exist blocks $B_i$ and $B_j$ such that $u, u', v, v' \in S'$ and $u, u' \in V(B_i)$ and $v, v' \in V(B_j)$, with a chain of blocks between $B_i$ and $B_j$ (denoted as $B_{i+1},B_{i+2},\dots, B_{j-1}$) such that $u_k\in S'\cap V(B_k)$ for $k=i+1,i+2,\dots,j-1$.

    Hence $\cup_{h=i}^{j} V(B_h)\subset \langle S'\setminus\{w\}\rangle$, or simply we can say that  $w\in \langle S'\setminus\{w\}\rangle$ for $w\in\{u,u',v,v'\}$. This implies that 
$S'$ is H-dependent.  
 In both cases, we conclude that $S'$ cannot be H-independent, contradicting the assumption that  $|S'| \geq \ell + 2$. Therefore, the largest H-independent set in $G$ has size $\ell + 1$, proving that $ h_\Delta(G) = \ell + 1 $.

\item[(b)] Let $G$ be a chordal graph containing blocks $B_1,B_2,B_3,\ldots, B_\ell$ and $B'_1,B'_2,B'_3,\ldots, B'_k$, where the blocks $B_i$ are complete for $i=1,2,\dots,\ell$ and the blocks $B'_j$ are non-complete for $j=1,2,\ldots,k$. Let  $G'$ be the induced subgraph of $G$ formed by $\bigcup_{i=i}^{k} V(B'_j)$. 

Let $S\subseteq V(G)$ be the set $S=\{$maximum independent vertex set of $G'\}\cup\{u_i\in V(B_i)\mid i=1,2,\ldots, \ell\}$. Then $|S| =\alpha(G')+\ell $. We can easily see that $a\notin\langle S\setminus \{a\}\rangle$ for any $a\in S$, and $V(B_i)\setminus S\nsubseteq\langle S\setminus \{b\}\rangle$ for any $b\in V(B_i)\cap S$. Therefore, $S$ is H-independent.

Next, we need to prove that there is no H-independent set of size greater than $\alpha(G')+\ell $. Suppose,  for the sake of contradiction, that there exists an H-independent set $S'$ with $|S'|>\alpha(G')+\ell$.
Consider the following cases:

             	\textbf{Case 1:} $S'$ contains more than $\alpha(G')$ vertices from the non-complete blocks. 
Then, $S'$ must contain at least two vertices each from at least two distinct non-complete blocks, implying that a pair of adjacent vertices exists in $S'$. Let $B'_j$ be a non-complete block such that $ S' \cap V(B'_j)$ contains adjacent vertices.  
\begin{itemize}
	\item  If $S' \cap V(B'_j)$ contains at least three vertices, say $u,v,w\in V(B'_j)$, with either $uv\in E(G)$ or $vw\in E(G)$, assume $uv\in E(G)$. Then, $V(B'_j)\subseteq \langle u,v\rangle$, implying $w\in \langle S'\setminus\{w\}\rangle$. Thus, $ w \in \bigcap_{a\in S'}\langle S'\setminus\{a\}\rangle$, and consequently, $S'$ is H-dependent.

	\item   If  $S' \cap V(B'_j) = \{u, v\}$ with  $uv \in E(G)$, then  $V(B'_j) \subseteq \langle u, v \rangle$. Let  $B'_{j+1}, B'_{j+2}, \ldots$ be a sequence of non-complete blocks where $B'_{j}$ is adjacent to $B'_{j+1}$,  $B'_{j+1}$ is adjacent to $B'_{j+2}$, and so on. If there exists  $w',w'' \in V(B'_k) \cap S'$ such that $w_{k}w' \in E(G)$, where $w_k \in V(B'_{k-1}) \cap V(B'_{k})$ for some $k=j+1,j+2,\ldots$,  then $V(B'_j) \cup V(B'_{j+1})\cdots \cup V(B'_{k})\subseteq \langle  (V(B'_j) \cup V(B'_{j+1})\cdots \cup V(B'_{k-1})\cup{w'})\cap S'\rangle$.  
	Consequently, $w'' \in \langle S' \setminus \{w''\} \rangle$, implying  $w'' \in \bigcap_{a \in S'} \langle S' \setminus \{a\} \rangle$, and $S'$ is H-dependent.
If there is no $w_i' \in V(B'_i) \cap S'$ with $w_iw_i' \in E(G)$ and $w_i \in V(B'_{i-1}) \cap V(B'_{i})$ for $i=j+1,\ldots,k$, then the blocks $B'_j,B_{j+1},\ldots, B'_i$ contribute additional vertices to the independence set, which contradicts the assumption that $ |S'| > \alpha(G')$. Thus, there must exist another non-complete block with adjacent vertices in $S'$, and the argument above applies to that block.
If there is no neighbouring non-complete block or there exists a non-complete block $B'_i$ such that $V(B'_i) \cap S' = \emptyset$, for $i=j+1,\ldots,k$ then, since  $|S'| > \alpha(G')$, there must exist another non-complete block containing adjacent vertices in $S'$, and the process described above applies similarly to that block. 
\end{itemize}
Thus, in all cases,  $S'$ is  H-dependent.\\

       \textbf{Case 2:} $S'$ contains more than $\ell$ vertices from the complete blocks.  If $S'$ contains three vertices from a complete block, then clearly $S'$ is H-dependent.  So $S'$ contains at most two vertices from each complete block. Let $u,u'\in V(B_i)\cap S'$, and let $B'_j$ be the nearest non-complete block (choose the pair with the shortest possible distance between them). Let $B_{i+1},B_{i+2},\ldots, B_k$ be the blocks between $B_i$ and $B'_j$.
    
           \textbf{Subcase (i):}If $S'$contains one vertex from each of the blocks $B_{i+1},B_{i+2},\ldots, B_k$ (if it is a cut vertex, count it for a single block). If there is $w\in V(B'_j)\cap V(B_k)\cap S'$, then $w\in \langle S'\setminus\{w\} \rangle$ implies that  $S'$ is H-dependent. If $V(B'_j)\cap V(B_k)\cap S'=\emptyset$, then let $w'\in V(B'_j)\cap V(B_k)$. Clearly, $w'x\in E(G)$ for some $x\in V(B'_j)\cap S'$; otherwise, $B'_j$ would not be a chordal component. Then $V(B'_j)\subseteq \langle S'\setminus\{y\}\rangle$ for $y\in V(B'_j)\cap S'$ and $y\neq x$ (such a $y$ exists since $| V(B'_j)\cap S'|\geq 2 $). This implies $S'$ is H-dependent.\\
           \textbf{Subcase (ii):}If $S'$ does not contains vertices from  some of the blocks among $B_{i+1},B_{i+2},\ldots, B_k$. That is, there is a subsequence of blocks $B_i, B_{i+1},\ldots, B_i'$ such that $S'$ contains two vertices from blocks $B_i$ and $B_i'$, and one vertex each from the blocks between $B_i$ and $B_i'$. Then clearly $S'$ is H-dependent since for $u\in V(B_i)\cap S'$, $u\in \langle S'\setminus \{u\}\rangle$.

 In both cases, we conclude that $S'$ cannot be H-independent, contradicting the assumption that  $|S'| >\alpha(G')+\ell$. Therefore, the largest H-independent set in $G$ has size $\alpha(G')+\ell$, proving that $ h_\Delta(G) = \alpha(G')+\ell $. 
\end{enumerate}		
  Hence the theorem.  \qed   
 \end{proof}

      \section{Radon number of the $\Delta-$ convexity}\label{section_Radon}

In this section, we investigate the Radon number within the framework of $\Delta$-convexity. We establish an upper bound for the Radon number that depends upon the number of triangles present in the graph and shows that this bound is the least upper bound. Furthermore, we determine the Radon number specifically for chordal graphs and note that it coincides with the Helly number for chordal graphs.

We start with the following observations, which arise directly from the definitions of $\Delta$-convexity and R-independent sets:

\begin{observation}
	If $G$ is a triangle-free graph, then the Radon number, $r_\Delta(G)=|V(G)|$.
\end{observation}
\begin{observation}
	If $G$ is a complete graph $K_n$ for $n>2$, then $r_\Delta(G)=2$.
\end{observation}
The following property follows directly from the definition of an R-dependent set: 
\begin{remark}\label{remark-Radon dependent} 
	If  $S$ is a set and $w \in S$ be such that $ w \in \langle S \setminus \{w\} \rangle $, then  $S$ is R-dependent.
\end{remark}
Next, we present some simple properties of the R-independent
set of a general graph, which will be utilized in subsequent results.
\begin{lemma}\label{lemma_Radon_characterization}
	Let $G$ be a graph. Then:
	\begin{enumerate}
		\item[(a)] The Radon number, $r_\Delta(G)$ satisfies $r_\Delta(G)\geq \alpha(G)$.
		\item[(b)] If $S$ is an R-independent set in $G$, then no three vertices of $S$ forms a $K_3$ in $G$.
		\item[(c)]If every pair of adjacent vertices forms a hull set in $G$, then $r_\Delta(G)=\max\{2,\alpha(G)\}$.
	\end{enumerate}
 \end{lemma}
	\begin{proof}
		\begin{enumerate}
			\item[(a)] This follows directly from Levi's inequality and Lemma~\ref{lemma_helly_characterization}.
   \item[(b)] Let $S$ be an R-independent set of $G$. Suppose $\{a,b,c\} \subseteq S$ forms an induced $K_3$ in $G$. Then $a \in \langle S \setminus \{a\} \rangle$ and the Remark~\ref{remark-Radon dependent}, implies that $S$ is R-dependent, a contradiction. Therefore, no three vertices of $S$ can form a $K_3$ in $G$.
   \item[(c)]	We proceed by considering the general case where  $G$ may be either a complete or a non-complete graph.
	
	First, consider the case where $G$ is a complete graph. In this case, the Radon number, $h_\Delta(G) = 2$. Since $\alpha(G) = 1$ for a complete graph, we observe that $\max\{2, \alpha(G)\} = 2$, so the lemma holds in this case.
	
	Now, suppose $G$ is not a complete graph. Let $S$ be the maximum independent vertex set of $G$. Thus $|S|=\alpha(G)$. Since every independent vertex set is R-independent, it follows that $S$ is R-independent. Now, consider any $S'\subseteq V(G)$ with $|S'|> \alpha(G)$. Then there must exist $u,v\in S'$ such that $uv\in E(G)$. By assumption, 
 $\langle u,v\rangle=V(G)$. 
 
 Now consider the partition of $S'=S_1'\cup S_2'$ where $S_1'=\{u,v\}$ and $S_2'=S\setminus S_1'$. Since $\langle S_1'\rangle\cap \langle S_2'\rangle=\langle S_2'\rangle$, this implies that $S'$ is R-dependent. Hence, the Radon number, $h_\Delta(G)=\alpha(G)$.\qed
	\end{enumerate}
	
	\end{proof}

From Lemmas \ref{chordal-convex} and \ref{lemma_Radon_characterization},  we immediately obtain the following corollary:
\begin{corollary}
    If $G$ is a 2-connected chordal graph, then $r_\Delta(G)=\max\{2,\alpha(G)\}$.
\end{corollary}
The following theorem provides an upper bound on the Randon number of a graph in $\Delta$-convexity:

\begin{theorem}
Let $G$ be a graph with $k$ triangles, and let $m$ denote the number of vertices not lying on any triangles. Then the Radon number, $r_\Delta(G)$ satisfies $r_\Delta(G)\leq  m + 2k$. 
\end{theorem}
\begin{proof}
	Let $S \subseteq V(G)$ be an  R-independent set in $G$. By Lemma \ref{lemma_Radon_characterization}, no three vertices in $S$ can form a triangle in $G$. Therefore, $S$ can contain at most two vertices from any triangle in $G$. Since $G$ contains $k$ triangles, the maximum number of vertices in $S$ that can be chosen from the triangles is $2k$.
	
	Next, consider the $m$ vertices of $G$ that do not lie on any triangles. Since these vertices are not part of any triangle, they are R-independent by default. Therefore, all of these $m$ vertices can be included in $S$. Therefore, the size of $S$ is bounded by $m + 2k$, implying that $r_\Delta(G) \leq m + 2k$.\qed
\end{proof}
\begin{remark}
    There exists a graph $G$ with $k$ triangles such that the Radon number, $r_\Delta(G) = m + 2k$, where $m$ denotes the number of vertices not lying on any triangle. For example, in Figure \ref{fig_rank}, consider the set $S = \{a_1, a_2, \ldots, a_k, b_1, b_2, \dots, b_k, d_1, d_2, \ldots, d_m\}$, where $|S| = m + 2k$. Then the following holds:
   $$ \left\langle S \setminus \{a_i\} \right\rangle = V(G) \setminus \{a_i, c_i\}, \quad \left\langle S \setminus \{b_i\} \right\rangle = V(G) \setminus \{b_i, c_i\}, \quad \text{and} \quad \left\langle S \setminus \{d_j\} \right\rangle = V(G) \setminus \{d_j\},$$

    for $i = 1, 2, \ldots, k$ and $j = 1, 2, \ldots, m$. For any partition $\{S_1, S_2\}$ of $S$, we have $\left\langle S_1 \right\rangle \cap \left\langle S_2 \right\rangle = \emptyset$, meaning $S$ is Radon independent.  Hence by Theorem 4, $r_{\Delta}(G)=m+2k$ follows.
\end{remark}
Furthermore, the Radon number can be made arbitrarily large. Specifically, for any natural number $n$, it is possible to construct a graph $G$ such that $r_\Delta(G) = n$, as stated in the following proposition.
\begin{proposition}\label{prop:rdelta}
    For any natural number $n > 1$, there exists a graph $G$ with $r_\Delta(G) = n$.
\end{proposition}

  \begin{proof}
    For $n = 2$, consider any complete graph $G$. In this case, $r_\Delta(G) = 2$.

    For $n \geq 3$, consider the graph $G$ depicted in Figure \ref{fig_delta}. Let $S = \{a_1, a_2, a_3, \ldots, a_n\}$. Then, $\langle S \rangle = V(G)$, $\langle S \setminus \{a_1\} \rangle = \{a_2, a_3, a_4, \ldots, a_n\}$, $\langle S \setminus \{a_2\} \rangle = \{a_1, a_3, a_4, \ldots, a_n\}$, and $\langle S \setminus \{a_i\} \rangle = \{a_1, a_2, b_1, a_3, b_2, a_4, \ldots, a_{i-1}, b_{i-2}, a_{i+1}, a_{i+2}, \ldots, a_n\}$ for $i = 3, 4, 5, \ldots, n$.

    For any partition $\{S_1, S_2\}$ of $S$, we have $\langle S_1 \rangle \cap \langle S_2 \rangle = \emptyset$, means that $S$ is R-independent.

    Next, we need to prove that there is no R-independent set of size greater than $n$. Suppose there exists an R-independent set $S'$ with $|S'| \geq n + 1$. It is clear that $S'$ does not contains all the vertices of a triangle.

    Now, consider the possible cases:
    \begin{itemize}
        \item[(a)] Suppose $S'$ contains one vertex from each triangle. Since $G$ contains only $n-1$ triangles, $S'$ must contain two vertices from at least two triangles, say $T_i$ and $T_j$. Let $a_i, b_i, a_j, b_j \in S'$, where $a_i, b_i \in T_i$ and $a_j, b_j \in T_j$. Then, clearly, $\langle S' \setminus a_i \rangle = \langle S' \setminus b_i \rangle = \langle S' \setminus a_j \rangle = \langle S' \setminus b_j \rangle = V(G)$. By Remark \ref{remark-Radon dependent}, $S'$ is R-dependent.
        
        \item[(b)] Suppose $S'$ does not contains vertices from some triangles. In this case, $S'$ contains two vertices from at least three triangles. Therefore, there exists a sequence of triangles $T_i, T_{i+1}, \ldots, T_{i+r}$ (where $i \in \{1, 2, \dots, k\}$ and $i + r \leq k$) such that $S'$ contains at least one vertex from each of these triangles and two vertices from at least two triangles, say $T_x$ and $T_y$ (with $x, y \in \{i, i+1, \dots, i+r\}$). Let $a_x, b_x \in T_x$ and $a_y, b_y \in T_y$. Then $a_x \in \langle S' \setminus a_x \rangle$, $b_x \in \langle S' \setminus b_x \rangle$, $a_y \in \langle S' \setminus a_y \rangle$, and $b_y \in \langle S' \setminus b_y \rangle$. By Remark \ref{remark-Radon dependent}, $S'$ is R-dependent.
    \end{itemize}
    
    Hence, for the graph in Figure \ref{fig_delta}, $r_\Delta(G) = n$.\qed
\end{proof}

  The following theorem provides the exact values for the Radon number of a chordal graph in $\Delta$-convexity.
\begin{theorem}
    Let $G$ be a chordal graph. Then: 
     \begin{enumerate}
        \item[(a)] If $G$  is  a block graphs with $\ell$ blocks, then $r_\Delta(G) = \ell + 1$. 
      \item[(b)]  If $G$  is  not a block graph, $i.e.$, $G$ contains non-complete blocks. Let $B_1,B_2,\dots, B_\ell$ be the complete blocks and $B'_1,B'_2,\dots, B'_k$ be the non-complete blocks. Let $G'$ be the induced subgraph of $G$ formed by $\bigcup_{i=i}^{k} V(B'_j)$. Then $r_\Delta(G) =\alpha(G')+\ell $.
    \end{enumerate}
\end{theorem}
 \begin{proof}
      \begin{enumerate}
          \item[(a)]     Let $G$ be a chordal graph where all the blocks are complete. That is, $G$ is a block graph. Let $B_1,B_2,\ldots,B_\ell$ be the blocks in $G$, where each block $B_i$ is complete for $i=\{1,2,\ldots,\ell\}$. Consider the set $S=\{u_1,u_1',u_2,\ldots,u_\ell\}$ such that $\{u_1,u_1'\}\in V(B_1)$ and $u_i \in V(B_i)$ for $i\in \{2,3,\dots,\ell\}$, and no three vertices of $S$ forms a $K_3$ in $G$.

       For any partition $\{S_1,S_2\}$ of $S$, if $u_1\in S_1$ and $u_1'\in S_2$, then $\langle S_1 \rangle=S_1$ and  $\langle S_2 \rangle=S_2 $, implying  $\langle S_1 \rangle \cap \langle S_2 \rangle=\emptyset$. Thus, $S$ is R-independent.

  Now, we need to show that there does not exist an R-independent set in $G$ that contains more than $\ell+1$ elements. Assume, to the contrary, that there exists an R-independent set $S'\subseteq V(G)$ such that $|S|\geq\ell+2$. Clearly, $S'$ cannot contain three vertices that form a $K_3$ in $G$. 
  
  Consider the following cases:
      
          \textbf{Case 1:} Suppose $S'$ contains at least one vertex from each block.
          
           Since $|S'|\geq \ell+2$, $S'$ contains two vertices from at least two blocks, say  $B_i$ and $B_j$ with $i, j \in \{1, 2, \dots, \ell\}$. Let $u, u' \in V(B_i)$ and $v, v' \in V(B_j)$, where $u, u', v, v' \in S'$. Then, $\langle S'\setminus\{u\}\rangle =\langle S'\setminus\{u'\}\rangle= \langle S'\setminus\{v\}\rangle =\langle S'\setminus\{v'\}\rangle = V(G)$. By Remark \ref{remark-Radon dependent}, $S'$ is $R$-dependent.\\
             \textbf{Case 2:} $S'$ does not contains vertices from some blocks in $G$. 

   Since $|S'| \geq \ell + 2$ and $S'$ does not contain more than two vertices from any single block, let $p$ be the number of blocks from which exactly two vertices are in $S'$, and $q$ be the number of blocks from which no vertices are in $S'$. Clearly, $p > q$.
   Thus, there exist blocks $B_i$ and $B_j$ such that $u, u', v, v' \in S'$ and $u, u' \in V(B_i)$ and $v, v' \in V(B_j)$, and a chain of blocks between $B_i$ and $B_j$ (denoted as $B_{i+1},B_{i+2},\ldots, B_{j-1}$) such that $u_k\in S'\cap V(B_k)$ for $k=i+1,i+2,\ldots,j-1$.

    Hence $\cup_{h=i}^{j} V(B_h)\subset \langle S'\setminus\{w\}\rangle$, or simply, $w\in \langle S'\setminus\{w\}\rangle$ for $w\in\{u,u',v,v'\}$. By remark ~\ref{remark-Radon dependent},
$S'$ is R-dependent.

 In both cases, we conclude that $S'$ cannot be R-independent, contradicting the assumption that  $|S'| \geq \ell + 2$. Therefore, the largest R-independent set in $G$ has size $\ell + 1$, proving that $ r_\Delta(G) = \ell + 1 $.
 
 \item[(b)] Let $G$ be a chordal graph containing blocks $B_1,B_2,\dots, B_\ell$ and $B'_1,B'_2,\dots, B'_k$, where the blocks $B_i$ are complete for $i=1,2,\dots,\ell$ and the blocks $B'_j$ are non-complete for $j=1,2,\dots,k$. Let  $G'$ be the induced subgraph of $G$ formed by $\bigcup_{i=i}^{k} V(B'_j)$. 

Let $S\subseteq V(G)$ be the set $S=\{$maximum independent vertex set of $G'\}\cup\{u_i\in V(B_i)\mid i=1,2,\ldots, \ell\}$. Then $|S| =\alpha(G')+\ell $. 

Next, we need to prove that there is no R-independent set of size greater than $\alpha(G')+\ell $. Suppose,  for the sake of contradiction, that there exists an R-independent set $S'$ with $|S'|>\alpha(G')+\ell$.
Consider the following cases:

	 \textbf{Case 1:} $S'$ contains more than $\alpha(G')$ vertices from the non-complete blocks. That is, $S'$ contains at least three vertices from some $B'_j$, say $u,v,w\in V(B'_j)$, with either $uv\in E(G)$ or $vw\in E(G)$. Suppose $uv\in E(G)$, then $\langle u,v\rangle=V(B'_i)$. So $w\in \langle S'\setminus\{w\}\rangle$, now we choose the partition of $S'$ as $S_1'=\{u,v\}$ and $S_2'=S'\setminus\{u,v\}$ which implies that $ w \in\langle S_1'\rangle\cap \langle S_2'\rangle$. Hence, $S'$ is R-dependent.\\
       \textbf{Case 2:} $S'$ contains more than $\ell$ vertices from the complete blocks.  If $S'$ contains three vertices from a complete block, then clearly $S'$ is R-dependent.  So $S'$ contains at most two vertices from each complete block. Let $u,u'\in V(B_i)\cap S'$, and let $B'_j$ be the nearest non-complete block to $B_i$. Also, Let $B_{i+1},B_{i+2},\ldots, B_k$ be the blocks between $B_i$ and $B'_j$ (shortest chain of blocks).\\
      
           \textbf{Subcase (i):} If $S'$contains one vertex from each of the blocks $B_{i+1}, B_{i+2},\ldots, B_k$ (if it is a cut vertex, count it for a single block). If there is $w\in V(B'_j)\cap V(B_k)\cap S'$, then $w\in \langle S'\setminus\{w\} \rangle$ implies that  $S'$ is R-dependent. If $V(B'_j)\cap V(B_k)\cap S'=\emptyset$, then let $w'\in V(B'_j)\cap V(B_k)$. Clearly, $w'x\in E(G)$ for some $x\in V(B'_j)\cap S'$; otherwise, $B'_j$ would not be a chordal component. Then $V(B'_j)\subseteq \langle S'\setminus\{y\}\rangle$ for $y\in V(B'_j)\cap S'$ and $y\neq x$ (such a $y$ exists since $| V(B'_j)\cap S'|\geq 2 $). That is, $y\in \langle S'\setminus\{y\}\rangle$. This implies $S'$ is 
           R-dependent.\\
           
           \textbf{Subcase (ii):}If $S'$ does not contains vertices from  some of the blocks among $B_{i+1},B_{i+2},\ldots, B_k$. That is, there is a subsequence of blocks $B_i, B_{i+1},\ldots, B_i'$ such that $S'$ contains two vertices from blocks $B_i$ and $B_i'$, and one vertex each from the blocks between $B_i$ and $B_i'$. Then clearly $S'$ is R-dependent, since for $u\in V(B_i)\cap S'$, $u\in \langle S'\setminus \{u\}\rangle$.
      
       In both cases, we conclude that $S'$ cannot be R-independent, contradicting the assumption that  $|S'| >\alpha(G')+\ell$. Therefore, the largest R-independent set in $G$ has size $\alpha(G')+\ell$, proving that $ r_\Delta(G) = \alpha(G')+\ell $. \qed
       \end{enumerate}
\end{proof}
 \section{Rank of the $\Delta$-convexity}\label{section_rank}
In this section, we delve into the concept of rank within the framework of $\Delta$-convexity. Similar to the previous sections, we establish an upper bound for the rank based on the number of triangles in the graph and demonstrate that this bound is the least possible upper bound.  Moreover, we determine the rank of block graphs. Interestingly, for block graphs, we observe that the rank is identical to both the Helly number and the Radon number in the context of $\Delta$-convexity.

 \begin{observation}
     If $G$ is a triangle free graph, then the rank, $d_\Delta(G)=|V(G)|$
 \end{observation}
 \begin{observation}
	If $G$ is a complete graph $K_n$ for $n>2$, then $d_\Delta(G)=2$.
\end{observation}
The following property follows directly from the definition of a convexly dependent set: 
\begin{remark}\label{remark-rank dependent} 
	If  $S$ is a set and $w \in S$ such that $ w \in \langle S \setminus \{w\} \rangle $, then  $S$ is convexly dependent.
\end{remark}
Next, we present some simple properties of the convexly independent set of a general graph, which will be used in subsequent results.
\begin{lemma}\label{lemma_Rank_Property}
For a graph $G$, the  following properties holds. 
    \begin{enumerate}
        \item [(i)] The rank, $d_\Delta(G)$ satisfies $d_\Delta(G) \geq\alpha(G)$.
        \item [(ii)] If $S$ is a convexly independent set, then no three vertices of $S$ forms a $K_3$ in $G$.
        \item [(iii)] If every two adjacent vertices of a graph $G$ form a hull set of $G$, then $d_\Delta(G)=\max\{2,\alpha(G)\}$.
    \end{enumerate}
    \end{lemma}
    \begin{proof}
        \begin{enumerate}
            \item [(i)]  Let $S$ be a maximum independent vertex set of the graph $G$, with $|S| = \alpha(G)$. Clearly, $\langle S \rangle = S$ and $ a \notin\langle S \setminus \{a\} \rangle$ for each $a \in S $.  This implies that $S$ is convexly independent. Therefore, $d_\Delta(G) \geq \alpha(G)$.
            
            \item [(ii)]Suppose $S \subseteq V(G)$ is a convexly independent set of $G$. Then by definition, $a \notin \langle S \setminus \{a\} \rangle$ for each $a \in S$. Let $S=
            \{a_1,a_2,a_3,\ldots ,a_n\};n \geq 3$ such that $S$ contains three vertices which forms a triangle $K_3$ in $G$. Then the removal of each vertex $a_i \in K_3$ is contained in $ \langle S \setminus \{a_i\} \rangle$ , contradicting our assumption that $S$ is convexly independent. Hence no three vertices of $S$ form a $K_3$ in $G$.
            
            \item [(iii)] Let $G$ be a graph such that for every $uv \in V(G),\langle u,v \rangle=V(G).$  Suppose $G$ is complete. Consider a subset $S \subseteq V(G)$ of cardinality 2. Let $S=\{u_1,u_2\}$. Clearly $u_1 \notin \langle S \setminus \{u_1\} \rangle$, $u_2 \notin \langle S \setminus \{u_2\} \rangle$ and $S$ is convexly independent. If possible, let $|S| \geq 3.$ Suppose $S=\{u_1,u_2,u_3\}$. If we consider $\langle S \setminus \{u_1\} \rangle$, the remaining vertices of $S$ can generate $V(G)$ and definitely $u_1 \in \langle S \setminus \{u_1\} \rangle$ implies that $S=\{u_1,u_2,u_3\}$ is convexly dependent. Hence $|S|=2.$\\
            Suppose $G$ is a non complete graph such that for any $uv \in E(G), \langle u,v \rangle =V(G).$ A maximum independent set $S$ is convexly independent and $d_{\Delta}(G)\geq \alpha (G)$ holds. If $|S|>\alpha(G)\geq 2$, then there exists different $x,y\in S$ such that $xy\in E(G)$. For $a\in S$, $a\notin\{x,y\}$, we have $a\in \left\langle S\setminus\{a\}\right\rangle$ and $d_{\Delta}(G)=\alpha (G)$.
        \end{enumerate}
    \end{proof}
From Lemmas \ref{chordal-convex} and \ref{lemma_Rank_Property},  we immediately obtain the following:

\begin{corollary}
   If $G$ is a 2-connected chordal graph, then $d_\Delta(G)=\max\{2,\alpha(G)\}$.
\end{corollary}
\begin{theorem}
Let $G$ be a graph with $k$ triangles, and let $m$ denote the number of vertices not lying on any triangles. Then $d_\Delta(G)\leq  m + 2k$. 
\end{theorem}
\begin{proof}
	A convexly independent set  $S\subseteq V(G)$ cannot contain more than two vertices from any triangle in $G$, limiting the number of vertices in $S$ that can be chosen from the $k$ triangles to at most $2k$.
    
     The remaining $m$ vertices, which do not lie on any triangles, are also convexly independent. Therefore, $S$ can include all of these $m$ vertices. Consequently, $d_\Delta(G) \leq m + 2k$.\qed
\end{proof}
\begin{remark}
    There exists a graph $G$ with $k$ triangles such that the rank, $d_\Delta(G) = m + 2k$; where $m$ is the number of vertices not lying on any triangle. For example, in Figure \ref{fig_rank}, consider the set $S = \{a_1, a_2, \ldots, a_k, b_1, b_2, \ldots, b_k, d_1, d_2, \ldots, d_m\}$, where $|S| = m + 2k$. Then the following holds:
   $$ \left\langle S \setminus \{a_i\} \right\rangle = V(G) \setminus \{a_i, c_i\}, \quad \left\langle S \setminus \{b_i\} \right\rangle = V(G) \setminus \{b_i, c_i\}, \quad \text{and} \quad \left\langle S \setminus \{d_j\} \right\rangle = V(G) \setminus \{d_j\}.$$ So $S$ is convexly independent.

   If any subset $S' \subseteq V(G)$ has  $|S'| > m + 2k$, it must include all vertices of at least one triangle, making $S'$ convexly dependent. Thus, $d_\Delta(G) = m + 2k$.
\end{remark}
Moreover, the rank can be made arbitrarily large. For any natural number $n$, there exists a graph $G$ with  $d_\Delta(G) = n$, as stated in the following proposition.
    \begin{proposition}\label{prop:rankdelta}
		For any natural number $n>1$, there exists a graph $G$ with $d_\Delta(G)= n$.
	\end{proposition}

\begin{proof}
 For $n=2$, consider any complete graph $G$. In this case, $d_\Delta(G)=2$.
		
		For $n\geq 3$, consider the graph $G$ shown in Figure \ref{fig_delta}. Let $S=\{a_1,a_2,\dots,a_n\}$, with $S$ convexly independent because $a_i\notin \langle S\setminus\{a_i\}\rangle$ for each $i$. Thus $d_\Delta(G)\geq n$.
		 
		 Now, we need to prove that there is no convexly independent set of size greater than $n$. Suppose there exists a convexly independent set $S'$ with $|S'|\geq n+1$. It is clear that $S'$ does not contains all vertices of a triangle.
   
		Consider the possible cases:
		 \begin{itemize}
		 	\item[(a)]Suppose $S'$ contains one vertex from each triangle. Since $G$ contains only $n-1$ triangles, $S'$ must contain two vertices from at least two triangles, say $T_i$ and $T_j$. Let  $a_i,b_i,a_j,b_j \in S'$, where $a_i,b_i \in T_i$ and $a_j,b_j \in T_j$. Then $\langle S'\setminus{a_i}\rangle=\langle S'\setminus{b_i}\rangle=\langle S'\setminus{a_j}\rangle=\langle S'\setminus{b_j}\rangle=V(G)$. By Remark \ref{remark-rank dependent},  $S'$ is convexly dependent. 
		 	\item[(b)]Suppose $S'$ does not contains vertices from some triangles. In this case, $S'$ contains two vertices from at least three triangles. Therefore, there exists a sequence of triangles $T_i,T_{i+1},\ldots, T_{i+r}$ (with $i\in \{1,2,\ldots, k\}$ and $i+r\leq k$) such that $S'$ contains at least one vertex from each of these triangles and two vertices from at least two triangles, say $T_x$ and $T_y$, (with $x,y\in\{i,i+1,\ldots,i+r\}$). Let $a_x,b_x \in T_x$ and $a_y,b_y \in T_y$. Then $a_x\in\langle S'\setminus{a_x}\rangle$, $b_x\in\langle S'\setminus{b_x}\rangle$ , $a_y\in \langle S'\setminus{a_y}\rangle$ and $b_y\in \langle S'\setminus{b_y}\rangle$. This
      makes $S'$ convexly dependent. 
		 \end{itemize}
		Hence, for the graph in Figure~\ref{fig_delta},  $d_\Delta(G)=n$. \qed
\end{proof}

For block graphs, the rank depends on the number of blocks in the graph. We present the following results for block graphs:
\begin{theorem}
If $G$ is a block graph with blocks $B_1,B_2,\dots, B_\ell$, then $d_\Delta(G) = \ell+1$.
 \end{theorem} 
 \begin{proof}
 Let $G$ be a block graph with blocks $B_1,B_2,\dots,B_\ell$. Consider a set $S\subseteq V(G)$ where $S=\{u_1,u_1',u_2,\ldots,u_\ell\}$, with $\{u_1,u_1'\} \in V(B_1)$ and $u_i \in V(B_i)$ for $i=2,3,4,\ldots,\ell$, such that no three vertices of $S$ forms a $K_3$ in $G$. Then, for each $a \in S$, we have $a \notin \langle S \setminus \{a\}\rangle$, making  $S$ convexly independent.

       Next, we need to show that there does not exist a convexly independent set in $G$ that contains more than $\ell+1$ vertices. Assume, to the contrary, that there exists a convexly independent set $S'\subseteq V(G)$ such that $|S|\geq \ell+2$. By Lemma \ref{lemma_Rank_Property},  $S'$ includes two vertices from at least two blocks and not more than two vertices from any single block. \\
       Consider the following cases:

 \textbf{Case 1:} $S'$ contains at least one vertex from each block. 
 
 Since $|S'|\geq \ell+2$, $S'$ must contain two vertices from at least two blocks, say  $B_i$ and $B_j$, where $i, j \in \{1, 2,\dots, \ell\}$. Let $u, u' \in V(B_i)$ and $v, v' \in V(B_j)$, where $u, u', v, v' \in S'$. Then, $\langle S'\setminus\{u\}\rangle =\langle S'\setminus\{u'\}\rangle= \langle S'\setminus\{v\}\rangle =\langle S'\setminus\{v'\}\rangle = V(G)$, implying $S'$ is convexly dependent.\\

  \textbf{Case 2:} $S'$ does not contains vertices from some blocks of $G$. 

   Since $|S'| \geq \ell + 2$ and $S'$ does not contain more than two vertices from any single block, let $p$ be the number of blocks in $G$ from which exactly two vertices are in $S'$, and $q$ be the number of blocks in $G$ from which no vertices are in $S'$. Clearly, $p > q$.
   Thus, there exists blocks $B_i$ and $B_j$ such that $u, u', v, v' \in S'$ with $u, u' \in V(B_i)$ and $v, v' \in V(B_j)$, and a chain of blocks between $B_i$ and $B_j$ (denoted as $B_{i+1},B_{i+2},\ldots, B_{j-1}$) such that $u_k\in S'\cap V(B_k)$ for $k=i+1,i+2,\ldots,j-1$. Hence $\cup_{h=i}^{j} V(B_h)\subset \langle S'\setminus\{w\}\rangle$, or simply, $w\in \langle S'\setminus\{w\}\rangle$ for $w\in\{u,u',v,v'\}$, making
$S'$ convexly dependent. 
 In both cases, we conclude that $S'$ cannot be convexly independent, violating our assumption that  $|S'| \geq \ell + 2$. Therefore, the largest convexly independent set in $G$ has size $\ell + 1$, establishing that $ d_\Delta(G) = \ell + 1 $.
 \qed
\end{proof}

\section*{Concluding Remarks}

It may be noted that, for most of the convexity in discrete structures, especially in graphs, the Radon number is one more than the Helly number. In this paper, we have shown that the Helly and Radon numbers are the same for chordal graphs with respect to $\Delta$-convexity.  We also observe that for all the small-sized graphs that we have examined so far, the Helly and Radon numbers are equal. This observation makes the $\Delta$-convexity rather special compared to most of the other graph convexities. Hence, we pose the following as a conjecture.

\begin{conjecture}
    If $G$ is a non-trivial connected graph, then $h_\Delta(G) = r_\Delta(G)$.
\end{conjecture}
\subsubsection*{Acknowledgment:} Arun Anil acknowledges the financial support from the University of Kerala, for providing the University Post Doctoral Fellowship (Ac EVII 5911/2024/UOK dated 18/07/2024). 
\bibliographystyle{amsplain}
\bibliography{deltahull}

 \end{document}